\documentclass[11pt,reqno]{amsart}
\usepackage{amsmath,amssymb,graphicx,amscd}
\usepackage{mathrsfs}
\usepackage{mathrsfs}
\usepackage{dsfont}

\newcommand{\CC}{{\mathbb C}}

\renewcommand{\span}{{\mathrm{span}}}

 \DeclareMathOperator{\Prob}{\mathbb{P}}   %probability
\newcommand{\Id}{{\mathrm{Id}}}
\newcommand{\ii}{{\mathrm{i}}}
\newcommand{\law}{\overset{\mbox{\rm \scriptsize law}}{=}}

\newtheorem{thm}{Theorem}[section]
\newtheorem{cor}[thm]{Corollary}

\theoremstyle{definition}

\theoremstyle{remark}
\newtheorem*{rem}{Remark}

\numberwithin{equation}{section}

\begin{document}
\title[On the characteristic polynomial on compact groups]
{The characteristic polynomial on compact groups with Haar measure :\\ some equalities in law}

\author{P. Bourgade}
 \address{Laboratoire de Probabilit\'es et Mod\'eles Al\'eatoires \\
Universit\'e Pierre et Marie Curie, et C.N.R.S. UMR 7599 \\ 175,
rue du Chevaleret \\ F-75013 Paris, France}
 \email{bourgade@enst.fr}

\author{A. Nikeghbali}
 \address{Institut f\"ur Mathematik,
 Universit\"at Z\"urich, Winterthurerstrasse 190,
 CH-8057 Z\"urich,
 Switzerland}
 \email{ashkan.nikeghbali@math.unizh.ch}

\author{A. Rouault}
 \address{Universit\'e  Versailles-Saint Quentin, LMV,
 B\^atiment Fermat, 45 avenue des Etats-Unis,
78035 Versailles Cedex}
 \email{rouault@fermat.math.uvsq.fr}

\date{14 June 2007}

\maketitle

\begin{center}
Presented by Marc Yor.
\end{center}
\vspace{0.3cm}

\begin{center}
\begin{minipage}[c]{0.8\linewidth}
\footnotesize
{\sc Abstract.}
This note presents some equalities in law for $Z_N:=\det(\Id-G)$,
where $G$ is an element of a subgroup of the set of unitary matrices of size $N$, endowed with its unique probability Haar measure.
Indeed, under some general conditions,
$Z_N$ can be decomposed as a product of independent random variables, whose laws are explicitly known.
Our results can be obtained in two ways : either by a recursive decomposition of the Haar measure (Section 1) or by previous results
by Killip and Nenciu (\cite{KillipNenciu}) on orthogonal polynomials with respect to some measure on the unit circle
(Section 2). This latter method leads naturally to a study of determinants of a class of principal submatrices.
\normalsize
\end{minipage}
\end{center}
\vspace{0.5cm}

\begin{center}
\begin{minipage}[c]{0.8\linewidth}
\footnotesize
{\sc R\'esum\'e.} Cette note pr\'esente quelques \'egalit\'es en loi pour
$Z_N:=\det(\Id-G)$, o\`u $G$ est un sous-groupe de l'ensemble des matrices unitaires de taille $N$, muni
de son unique mesure de Haar normalis\'ee. En effet, sous des conditions assez g\'en\'erales,
$Z_N$ peut \^etre d\'ecompos\'e comme le produit de variables al\'eatoires ind\'ependantes, dont on connait la loi explicitement.
Notre r\'esultat peut \^etre obtenu de deux mani\`eres : soit par une d\'ecomposition r\'ecursive de la mesure de Haar (Partie 1)
soit en utilisant un r\'esultat de Killip et Nenciu (\cite{KillipNenciu}) \`a propos des polyn\^omes orthogonaux
relativement \`a une certaine mesure sur le cercle unit\'e (Partie 2). Cette derni\`ere m\'ethode nous conduit naturellement \`a
l'\'etude des d\'eterminants de certaines sous-matrices.
\normalsize
\end{minipage}
\end{center}
\vspace{1cm}

In this note, $\langle a,b\rangle$  denotes the Hermitian product of  two elements $a$ and $b$ in $\CC^N$ (the dimension is implicit).

\section{A recursive decomposition, consequences} \label{Section1}

\subsection{The general equality in law.} Let $\mathcal{G}$ be a subgroup of $U(N)$, the group of unitary matrices of size $N$.
Let $(e_1,\dots,e_N)$ be an orthonormal basis of $\CC^N$ and $\mathcal{H}:=\{H\in\mathcal{G}\mid H(e_1)=e_1\}$, the subgroup of $\mathcal{G}$
which stabilizes $e_1$. For a generic compact group $\mathcal{A}$, we write $\mu_{\mathcal{A}}$ for the unique Haar probability measure on $\mathcal{A}$.
Then we have the following Theorem.

\begin{thm}\label{thm:decomposition}
Let $M$ and $H$ be independent matrices, $M\in\mathcal{G}$ and $H\in\mathcal{H}$ with distribution $\mu_\mathcal{H}$.
Then $M H\sim\mu_\mathcal{G}$ if and only if $M(e_1)\sim f(\mu_{\mathcal{G}})$, where $f$ is the map $f : G\mapsto G(e_1)$.
\end{thm}

Let $\mathcal{M}$ be the set of elements of $\mathcal{G}$ which are reflections with respect to a hyperplane of $\CC^N$.
Define also
$$
g :
\left\{
\begin{array}{ccc}
\mathcal{H}&\to&U(N-1)\\
H&\mapsto&H_{\span(e_2,\dots,e_N)}
\end{array}
\right.,
$$
where $H_{\span(e_2,\dots,e_N)}$ is the restriction of $H$ to
$\span(e_2,\dots,e_N)$. Now suppose that  $\{G(e)\mid
G\in\mathcal{G}\}=\{M(e)\mid M\in\mathcal{M}\}$. Under this
additional condition the following Theorem can be proven, using
Theorem \ref{thm:decomposition} and elementary manipulations of
determinants.

\begin{thm}\label{thm:generaldeterminants}
Let $G\sim\mu_\mathcal{G}$, $G'\sim\mu_\mathcal{G}$ and $H\sim g(\mu_\mathcal{H})$ be independent.
Then $$\det(\Id_N-G)\law (1-\langle e_1,G'(e_1)\rangle)\det(\Id_{N-1}-H).$$
\end{thm}

\subsection{Examples : the unitary group, the group of permutations.}
Take $G=U(N)$. As all reflections with respect to a hyperplane of $\CC^N$ are elements of $G$,
one can apply Theorem \ref{thm:generaldeterminants}. The corresponding measures are the following.
\begin{enumerate}
\item The distribution $g(\mu_{\mathcal{H}})$ is clearly
$\mu_{U(N-1)}$.
\item $\langle e_1,G(e_1)\rangle$ is distributed as the first coordinate of a vector
of the $N-$dimensional unit complex sphere with uniform measure : $\langle e_1,G(e_1)\rangle\sim
e^{i\theta_N}\sqrt{\beta_{1,N-1}}$ with $\theta_n$ uniform on $(0,2\pi)$
and independent of $\beta_{1,N-1}$, a beta variable with  parameters 1 and $N-1$.
\end{enumerate}

Thus iterations of Theorem \ref{thm:generaldeterminants} lead to the following Corollary.

\begin{cor} (\cite{BHNY}) \label{unitarycase}
Let $G \in U(N)$ be $\mu_{U(N)}$ distributed. Then
$$
\det (\Id_N-G)\law \prod_{k=1}^N
\left(1-e^{\ii\theta_k}\sqrt{\beta_{1,k-1}}\right),
$$
with $\theta_1,\dots,\theta_N,\beta_{1,0},\dots,\beta_{1,N-1}$
independent random variables, the $\theta_k$'s uniformly distributed
on $(0,2\pi)$ and the $\beta_{1,j}$'s ($0\leq j\leq N-1$) being beta
distributed with parameters 1 and $j$ (by convention, $\beta_{1,0}$
is the Dirac distribution at 1).
\end{cor}

The group $\mathcal{S}_N$ of permutations of size $N$ gives another possible application.  Identify an element
$\sigma\in\mathcal{S}_N$ with the matrix $(\delta^j_{\sigma(i)})_{1\leq i,j\leq N}$ ($\delta$ is Kronecker's symbol).
As $\det(\Id_N-\sigma)$ is equal to 0, we prefer to deal with the group $\tilde {\mathcal{S}}_N$ of matrices
$(e^{\ii\theta_j}\delta^j_{\sigma(i)})_{1\leq i,j\leq N}$, with $\sigma\in\mathcal{S}_N$ and
$\theta_1,\dots,\theta_N$ independent uniform random variables on $(0,2\pi)$. Then the measures
corresponding to Theorem \ref{thm:generaldeterminants} are the following.
\begin{enumerate}
\item The distribution $g(\mu_{\tilde {\mathcal{S}}_N})$ is
$\mu_{\tilde{\mathcal{S}}_{N-1}}$.
\item $\langle e_1,G(e_1)\rangle$ is 0 with probability $1-1/N$ and $e^{\ii \theta}$ ($\theta$ uniform on $(0,2\pi)$) with probability $1/N$.
\end{enumerate}
As previously, iterations of Theorem \ref{thm:generaldeterminants} give the following result.

\begin{cor} \label{permutationscase}
Let $S_N \in \tilde {\mathcal{S}}_N$ be $\mu_{\tilde {\mathcal{S}}_N}$
distributed. Then
$$
\det (\Id_N-S_N)\law \prod_{k=1}^N
\left(1-e^{\ii\theta_k}X_k\right),
$$
with $\theta_1,\dots,\theta_N,X_1,\dots,X_N$
independent random variables, the $\theta_k$'s uniformly distributed
on $(0,2\pi)$ and the $X_k$'s Bernoulli variables : $\Prob(X_k=1)=1/k$, $\Prob(X_k=0)=1-1/k$.
\end{cor}

\begin{rem}
Let $k_\sigma$ be the number of cycles of a random permutation of size $N$, with respect to the (probability) Haar measure.
Corollary \ref{permutationscase} allows us to recover the following celebrated result about the law of $k_\sigma$ :
$$k_\sigma\law X_1+\dots+X_N,$$
with the previous notations. Indeed, if a permutation $\sigma\in\mathcal{S}_N$ has $k_\sigma$ cycles with lengths
$l_1,\dots,l_{k_\sigma}$ ($\sum_k l_k=N$), then it is easy to see that under the Haar measure
$$\det(x \Id-\tilde{\mathcal{S}}_N)\law \prod_{k=1}^{k_\sigma}(x^{l_k}-e^{\ii\alpha_k})$$
with the $\alpha_k$'s independent and uniform on $(0,2\pi)$. Using the previous relation and the result of Corollary \ref{permutationscase}
we get
$$
\prod_{k=1}^N
\left(1-e^{\ii\theta_k}X_k\right)\law
\prod_{k=1}^{k_\sigma}(1-e^{\ii\alpha_k}).
$$
The equality of the Mellin transforms of the modulus of the above members easily implies the expected result :
$k_\sigma\law X_1+\dots+X_N$. Our discussion on the permutation group is closely related to the so-called Chinese
restaurant process and the Feller decomposition of the symmetric group (see, e.g. \cite{ABT}).
\end{rem}

\section{Characteristic polynomials as orthogonal polynomials}

We now show how Corollary \ref{unitarycase} can be obtained as a consequence of a result by
Killip and Nenciu (\cite{KillipNenciu}).

\subsection{A result by Killip and Nenciu.}
Let $\mathbb D$ be the open unit disk $\{ z \in \mathbb C : |z| < 1\}$ and $\partial\mathbb D$ the unit circle.
Let  $(e_1, \dots, e_N)$ be the canonical basis of  $\mathbb C^N$.
If $G \in  U(N)$, and if $e_1$ is  cyclic for $G$, the spectral measure for the pair
$(G, e_1)$ is the unique   probability $\nu$ on $\partial\mathbb D$ such that, for every integer $k \geq 0$
\begin{equation}
\label{moment}\langle e_1 , G^k e_1\rangle = \int_{\partial\mathbb D} z^k d\nu(z).
\end{equation}
In fact, we have the expression
$$
\nu = \sum_{j=1}^N \pi_j \delta_{e^{i\zeta_j}}$$
where $(e^{i\zeta_j}, j=1, \dots N)$ are the eigenvalues of  $G$ and where
$\pi_j = |\langle e_1, \Pi e_j\rangle|^2$ with $\Pi$  a unitary matrix  diagonalizing $G$.

The  relation (\ref{moment}) allows to define an isometry from $\mathbb C^N$ equipped with the basis
$(e_1, Ge_1, \cdots , G^{N-1}e_1)$  into the subspace of $L^2(\partial\mathbb D ; d\nu)$ spanned by the family
$(1, z, \dots, z^N)$. The endomorphism
$G$ is then a representation of the multiplication by $z$.

From the linearly independent family of monomials  $\{1, z, z^2, \dots, z^{N-1}\}$ in $L^2(\partial \mathbb D, \nu)$,
we construct an orthogonal basis $\Phi_0, \dots, \Phi_{N-1}$ of monic polynomials by the Gram-Schmidt procedure.
The $N^{th}$ degree polynomial obtained this way is
\[\Phi_N (z) =\prod_{j=1}^N (z - e^{i\zeta_j}),\]
i.e. the characteristic polynomial of $G$. The $\Phi_k$'s  ($k= 0, \dots, N$)
obey the Szeg\"o recursion relation:
\begin{equation}
\label{Szego}
\Phi_{j+1}(z) = z\Phi_j(z) - \bar\alpha_j \Phi_j^*(z)
\end{equation}
where $\Phi_j^*(z) = z^j\!\ \overline{\Phi_j(\bar z^{-1})}$.
The coefficients $\alpha_j's$ ($j \geq 0)$ are called   Schur  or Verblunsky coefficients and satisfy the condition
$\alpha_0, \cdots , \alpha_{N-2} \in \mathbb D$ and
$\alpha_{N-1} \in \partial\mathbb D$.
There is a bijection between this set of coefficients and the set of spectral probability measures $\nu$ (Verblunsky's theorem).
If $G\sim\mu_{U(N)}$, then we know the exact distribution of the Verblunsky coefficients :

\begin{thm}\label{thm:ThKN}(Killip and Nenciu \cite{KillipNenciu})
\label{ThKN} Let $G \in U(N)$ be $\mu_{U(N)}$ distributed.
The Verblunsky parameters $\alpha_0, \cdots , \alpha_{N-2}, \alpha_{N-1}$
are independent and the density of $\alpha_j$ for $j \leq N-1$ is
$$\frac{N-j-1}{\pi} \left(1 -
|z|^2\right)^{N-j-2}\mathds{1}_{\mathbb D}(z)$$
(for $j=N-1$ by convention this is the uniform measure on the unit circle).
\end{thm}

\subsection{Recovering Corollary \ref{unitarycase}.} For $z=1$, Szeg\"o's recursion (\ref{Szego}) can be written
\begin{equation}\label{AtOne}
\Phi_{j+1}(1)=\Phi_j(1)-\overline{\alpha_j}\,\overline{\Phi_j(1)}.
\end{equation}
Under the Haar measure for $G$, as $\alpha_j$ is independent of $\Phi_j(1)$ and
its distribution is invariant by rotation, (\ref{AtOne}) easily yields
$$
\Phi_{j+1}(1)\law (1-\alpha_j)\Phi_j(1).
$$
In particular, for $j=N-1$ we get by induction
\begin{equation}\label{Atj}
\det(\Id-G)\law\prod_{k=0}^{N-1}(1-\alpha_j).
\end{equation}
From the density for $\alpha_j$ given in Theorem \ref{thm:ThKN}
one can see that this is exactly the same result as Corollary \ref{unitarycase}.

\begin{rem}
A similar result holds for $SO(2N)$, and can be shown using either the method of Section 1 or the one in Section 2, with the corresponding
result by Killip and Nenciu for the Verblunsky coefficients on the orthogonal group \cite{KillipNenciu}.
\end{rem}
\subsection{Extension.} We  now consider the whole sequence of polynomials $\Phi_j, j \leq N$ for $j \leq N$
as a sequence of characteristic polynomials. For this purpose, we apply the Gram-Schmidt procedure to
$1, z, z^{-1}, z^2,\dots,z^{p-1}, z^{1-p}, z^p$ if $N=2p$ and to
$1, z, z^{-1}, z^2,\dots,z^{p}, z^{-p}$ if $N=2p+1$ in $L^2(\partial\mathbb D); d\nu)$.
In the resulting basis, the mapping $f(z) \mapsto zf(z)$ is represented by a so-called CMV matrix
(\cite{KillipNenciu} Appendix B, \cite{simon2006cmf}) denoted by ${\mathcal C}_N(G)$.
It is five-diagonal and conjugate to $G$. For $1 \leq j \leq N$ let ${\mathcal C}_N^{(j)}(G)$
the principal submatrix of order $j$ of ${\mathcal C}_N(G)$. It is known
(see for instance Proposition 3.1 in \cite{simon2006cmf}) that
\[\Phi_j (z) = \det \Big(z\Id_j - {\mathcal C}_N^{(j)}(G)\Big)\,.\]
From the recursion (\ref{AtOne}) and looking at the invariance of conditional distributions,
we see that
\begin{equation}\label{sequence}
\left(\det \Big(\Id_j - {\mathcal C}_N^{(j)}(G)\Big)\right)_{1 \leq j \leq N}
= \left(\Phi_j(1)\right)_{1 \leq j \leq N}
\law  \left(\prod_{l=0}^j(1- \alpha_l)\right)_{0 \leq j \leq N-1}.
\end{equation}
It allows a study of the process
$\big(\log \Phi_{\lfloor Nt\rfloor}(1)\ , \ t \in [0,1]\big)$
as a triangular array of (complex) independent random variables.
For $t=1$ the asymptotic behavior is presented in \cite{BHNY}
(see (\ref{CLT1} below). It is remarkable that for $t<1$, we do not need any normalization for the CLT.
\begin{thm}\begin{enumerate}
\item
As $N \rightarrow \infty$
\begin{equation}\label{Donskerbis}
\big(\log \det \big(\Id_j - {\mathcal C}_{\lfloor N t\rfloor}^{(j)}(G)\big); \ t \in [0,1)\big) \Rightarrow \big(
{\bf B}_{-\frac{1}{2}\log (1-t)} ; \ t \in [0,1)\big)\,, \end{equation}
where ${\bf B}$ is a standard complex Brownian motion
 and $\Rightarrow$ stands for the weak
convergence of distributions in the set of c\`adl\`ag functions on  $[0, 1)$,  starting from $0$, endowed with
the Skorokhod topology.
\item As $N \rightarrow \infty$, \begin{equation}\label{CLT1}\frac{\log \det(\Id_N - G)}{\sqrt{2\log N}}
\Rightarrow {\mathcal N_1 + \ii\mathcal N_2 }\end{equation}
where $\mathcal N_1$ and  $\mathcal N_2$ are independent standard normal and  independent of
$\bf B$,
and $\Rightarrow$ stands for the weak convergence of distributions
in $\mathbb C$.
\end{enumerate}
\end{thm}
This theorem can be proved using the Mellin-Fourier transform of the
$1-\alpha_j$'s and independence. This method may also be used to prove large deviations. It is the topic of a companion paper.
These results occur in similar way for other random determinants (see \cite{rou}).

\renewcommand{\refname}{References}

\end{document}